\documentclass[11pt]{amsart}
\usepackage{amssymb,amsmath,amsthm}
\newtheorem{theorem}{Theorem}[section]

\theoremstyle{definition}

\theoremstyle{remark}

\numberwithin{equation}{section}

\newcommand{\abs}[1]{\lvert#1\rvert}

\newcommand\Etwid{\overset {\text{\lower 3pt\hbox{$\sim$}}}E}
\newcommand\Ftwid{\overset {\text{\lower 3pt\hbox{$\sim$}}}F}
\newcommand\Qtwid{\overset {\text{\lower 3pt\hbox{$\sim$}}}Q}

\DeclareSymbolFont{AMSb}{U}{msb}{m}{n}
\DeclareMathSymbol{\Z}{\mathalpha}{AMSb}{"5A}

\DeclareMathSymbol{\nmid}{\mathrel}{AMSb}{"2D}
\DeclareSymbolFont{AMSb}{U}{msb}{m}{n}
\DeclareMathSymbol{\C}{\mathalpha}{AMSb}{"43}
\DeclareMathSymbol{\F}{\mathalpha}{AMSb}{"46}
\DeclareMathSymbol{\N}{\mathalpha}{AMSb}{"4E}
\DeclareMathSymbol{\Q}{\mathalpha}{AMSb}{"51}
\DeclareMathSymbol{\R}{\mathalpha}{AMSb}{"52}
\DeclareMathSymbol{\Z}{\mathalpha}{AMSb}{"5A}

\allowdisplaybreaks

\begin{document}
\newcommand{\beqs}{\begin{equation*}}
\newcommand{\eeqs}{\end{equation*}}
\newcommand{\beq}{\begin{equation}}
\newcommand{\eeq}{\end{equation}}
\newcommand\mylabel[1]{\label{#1}}
\newcommand\eqn[1]{(\ref{eq:#1})}
\newcommand\thm[1]{\ref{thm:#1}}
\newcommand\lem[1]{\ref{lem:#1}}
\newcommand\propo[1]{\ref{propo:#1}}
\newcommand\corol[1]{\ref{cor:#1}}
\newcommand\sect[1]{\ref{sec:#1}}
\newcommand\gausspoly[2]{\begin{bmatrix} #1 \\ #2\end{bmatrix}_q}
\newcommand\mytwid[1]{\overset {\text{\lower 3pt\hbox{$\sim$}}}#1}
\newcommand\leg[2]{\genfrac{(}{)}{}{}{#1}{#2}} 
\newcommand\tpi{(\pi_1,\pi_2)}
\newcommand\bpi{{\mathbf\pi}}
\newcommand\HLbr{\mbox{\scriptsize{HL-birank}}(\bpi)}          
\newcommand\gHLmr{\mbox{\scriptsize{gHL-multirank}}(\vec{\pi})}          
\newcommand\mcrI{\mbox{\scriptsize{multicrank-I}}(\vec{\pi})}          
\newcommand\mcrII{\mbox{\scriptsize{multicrank-II}}(\vec{\pi})}          
\newcommand\DYbr{\mbox{\scriptsize{Dyson-birank}}(\bpi)}          
\newcommand\bcra{\mbox{\scriptsize{bicrank}}_1(\bpi)}          
\newcommand\bcrb{\mbox{\scriptsize{bicrank}}_2(\bpi)}          
\newcommand\HL{\mbox{\scriptsize{HL}}}          
\newcommand\gHL{\mbox{\scriptsize{gHL}}}          
\newcommand\DY{\mbox{\scriptsize{D}}}          
\newcommand\FC{\mbox{\scriptsize{5C}}}          
\newcommand\srank{\mbox{\scriptsize{rank}}}          
\newcommand\fccrank{\mbox{\scriptsize{$5$-core-crank}}}          
\newcommand\fcbr{\mbox{\scriptsize{$5$-core-birank}}(\bpi)}          
\newcommand\ord{\mbox{ord}\,}
\newcommand\ORD{\mbox{ORD}\,}

\title[Biranks for Partitions into $2$ Colors]{
Biranks for Partitions into $2$ Colors}

\author{F.G. Garvan}
\address{Department of Mathematics, University of Florida, Gainesville,
Florida 32611-8105}
\email{frank@math.ufl.edu}          
\thanks{The author was supported in part by NSA Grant H98230-09-1-0051.
The main results of this paper were presented as part of a talk ``The Rank
and Crank of Partitions --- Dedicated to the Memory of Richard P. Lewis''
given on Friday, March 14, 2008 at the Gainesville Partitions, $q$-Series and
Modular Forms Conference.
}



\date{September 26, 2009}  %

\dedicatory{This paper is dedicated to the memory of my friend 
Richard Lewis (1942 -- 2007)}


\begin{abstract}
In 2003, Hammond and Lewis defined a statistic on partitions into $2$ colors
which combinatorially explains certain well known partition congruences
mod $5$. We give two analogs of Hammond and Lewis's birank statistic.
One analog is in terms of Dyson's rank
and the second uses the $5$-core crank due to Garvan, Kim and Stanton.
We discuss Andrews's bicrank statistic and how it may be extended.
We also generalize the Hammond-Lewis birank to a multirank for
multipartitions and the Andrews bicrank to a multicrank for
extended multipartitions. These both give combinatorial interpretations for
multipartition congruences modulo all primes $t>3$.
\end{abstract}

\maketitle

\section{Introduction} \label{sec:intro}

Hammond and Lewis \cite{Ha-Le} found some elementary results for 
$2$-colored partitions mod $5$.
Let $E(q) = \prod_{n=1}^\infty(1-q^n)$, and
$$
\sum_{n=0}^\infty p_{-2}(n) q^n = \frac{1}{E(q)^2},
$$
which is the generating function for pairs of partitions $(\pi_1,\pi_2)$
(or $2$-colored partitions).  Throughout this paper we refer to
such pairs of partitions as {\it bipartitions}. It is not hard to show that
\begin{equation}
p_{-2}(5n+2) \equiv 
p_{-2}(5n+3) \equiv 
p_{-2}(5n+4) \equiv 0 \pmod{5}.
\mylabel{eq:2colcongs}
\end{equation}
Hammond and Lewis \cite{Ha-Le} found a crank for these congruences. By crank
we mean a statistic that divides the relevant partitions into equinumerous
classes.
They define
\begin{equation}
\mbox{birank}(\pi_1,\pi_2) = \#(\pi_1) - \#(\pi_2),
\mylabel{eq:birankdef}
\end{equation}
where $\#(\pi)$ denotes the number of parts in the partition $\pi$. 
They show that the residue of the birank mod $5$ 
divides the bipartitions of $n$ into $5$ equal classes
provided $n\equiv2$, $3$ or $4\pmod{5}$.
The proof is elementary.  It relies on Jacobi's triple product identity and the 
method of \cite{Ga88}, which  uses roots of unity. 
We have found two other analogs. 

\bigskip
\noindent 
\textbf{First Analog - The Dyson-birank}

\noindent
Dyson \cite{Dy44} defined the \textit{rank} of a partition as
the largest part minus the number of parts. We define
\begin{equation}
\mbox{Dyson-birank}(\pi_1,\pi_2) = \mbox{rank}(\pi_1)
+ 2\,\mbox{rank}(\pi_2).
\mylabel{eq:Dysonbirankdef}
\end{equation}
In Section \sect{Dbirank}, we show that the residue of the Dyson-birank mod $5$ 
divides the bipartitions of $n$ into $5$ equal classes
provided $n\equiv2$, or $4\pmod{5}$. Unfortunately the Dyson-birank
does not work if $n\equiv3\pmod{5}$. Nonetheless, for the other
residue classes this is a surprising and deep result because
of the nature of the rank generating function. The proof
depends on known results for the rank mod $5$ due to
Atkin and Swinnerton-Dyer \cite{At-Sw}.

\bigskip
\noindent 
\textbf{Second Analog - The $5$-core-birank}

\noindent
In \cite{Ga-Ki-St} new statistics were defined in terms of $t$-cores which 
gave new combinatorial interpretations of Ramanujan's partition
congruences mod $5$, $7$ and $11$. For example, for a partition $\pi$
the \textit{$5$-core-crank} is defined as 
\begin{equation}
\mbox{$5$-core-crank}(\pi) = r_1 + 2 r_2 - 2 r_3 - r_4,
\mylabel{eq:GKScrankdef}
\end{equation}
where $r_j$ is the number of cells labelled $j$ in the $5$-residue
diagram of $\pi$. Then in \cite{Ga-Ki-St} we proved combinatorially that
the residue of the $5$-core-crank divides the partitions of $5n+4$
into $5$ equal classes.
We define
\begin{equation}
\mbox{$5$-core-birank}(\pi_1,\pi_2) = \mbox{$5$-core-crank}(\pi_1)
+ 2\,(\mbox{$5$-core-crank}(\pi_2)).
\mylabel{eq:GKSbirankdef}
\end{equation}
In Section \sect{5cbirank}, we show that the $5$-core-birank divides the 
bipartitions of $n$ into 
$5$ equal classes
for $n\equiv2$, $3$ or $4\pmod{5}$. This is quite a surprising result.
The proof relies on the $5$-dissection 
of the $5$-core-crank generating function for
$5$-cores. 

The crank of a partition is defined to be the largest part if
it contains no ones and otherwise it is the difference between number of parts larger
than the number of ones, and the number of ones.  The crank gives a combinatorial
of Ramanujan's partition congruences mod $5$, $7$ and $11$ and solves a
problem of Dyson \cite{Dy44}, \cite[p.52]{Dy96}. See \cite{An-Ga88}. This crank
is different to the $5$-core crank given in \cite{Ga-Ki-St}. It is natural
to ask whether there is a crank analog of the birank. This question has been
answered in part by Andrews \cite{An08}. In Section \sect{bicrank},
we consider Andrews result and how it may be extended. In \cite{An08},
Andrews also considered congruences for more general multipartitions. In
Section \sect{ghlbirank}, we give multipartition analogs of the Hammond-Lewis birank which
explain these more general congruences. In Section \sect{mcranks}, we extend
Andrews bicrank to multicranks of what we call extended multipartitions, and
give alternative explanations of our multipartition congruences. In Section \sect{end},
we close with some further problems.

\subsection*{Notation}
For a partition $\pi$ we denote the sum of parts by $\abs{\pi}$.
We will use the standard $q$-notation.
$$
(z;q)_n=(z)_n=
\begin{cases}
\prod_{j=0}^{n-1}(1-zq^j), & n>0 \\
1,                         & n=0,
\end{cases}
$$
and
$$
(z;q)_\infty=(z)_\infty = \lim_{n\to\infty} (z;q)_n
=\prod_{n=1}^\infty (1-z q^{(n-1)}),
$$
where $\abs{q}<1$. We will also you the following notation
for Jacobi-type theta-products.
$$
J_{a,m}(q) := (q^a;q^m)_\infty (q^{m-a};q^m)_\infty (q^m;q^m)_\infty.
$$


\section{The Hammond-Lewis Birank} \label{sec:HLbirank}

For completeness we include some details of the Hammond-Lewis birank.
For a bipartition $\bpi=\tpi$ we denote the sum of parts by
\beq
\abs{\bpi} = \abs{\pi_1} + \abs{\pi_2}.
\mylabel{eq:abspi}
\eeq
We denote the Hammond-Lewis birank by
\begin{equation}
\mbox{HL-birank}(\bpi) = \#(\pi_1) - \#(\pi_2),
\mylabel{eq:HLbirankdef}
\end{equation}
where $\#(\pi)$ denotes the number of parts in the partition $\pi$.
The HL-birank generating function is
\beq
\sum_{\bpi=(\pi_1,\pi_2)} z^{\HLbr} q^{\abs{\bpi}} = \frac{1}{(zq;q)_\infty (z^{-1}q;q)_\infty}.
\mylabel{eq:HLgenfunc}
\eeq           
We let $N_{\HL}(m,t,n)$ denote the number of bipartitions $\bpi=\tpi$
with HL-birank congruent to $m\pmod{t}$. 
Suppose $\zeta$ is primitive $5$th root of unity.
By letting $z=\zeta$ in \eqn{HLgenfunc}
and
using Jacobi's triple product identity, Hammond and Lewis found that
\begin{align}
\sum_{n=0}^\infty 
\sum_{k=0}^4 \zeta^k N_{\HL}(k,5,n) q^n
&= \frac{1}{(\zeta q;q)_\infty (\zeta^{-1}q;q)_\infty}
= \frac{(\zeta^2q;q)_\infty (\zeta^{-2}q;q)_\infty(q;q)_\infty}
        {(q^5;q^5)_\infty},\mylabel{eq:birankmod5}\\
&=(q^{25};q^{25})_\infty\left(
      \frac{1}{J_{1,5}(q^5)} + (\zeta + \zeta^{-1})
      \frac{q}{J_{2,5}(q^5)}
       \right).                      
\nonumber
\end{align}
Since the coefficient of $q^n$ on the right side of \eqn{birankmod5}
is zero when $n\equiv2$, $3$ or $4\pmod{5}$, Hammond and Lewis's main result
follows.
\begin{theorem}
\cite{Ha-Le}
\label{thm:thm1}
The residue of the HL-birank mod $5$
divides the bipartitions of $n$ into $5$ equal classes
provided $n\equiv2$, $3$ or $4\pmod{5}$.
\end{theorem}

We illustrate Theorem \thm{thm1} for the case $n=3$.
\beqs
\begin{array}{rr}
\mbox{Bipartitions of $3$} & \mbox{HL-birank (mod $5$)}\\
(3,-) & 1-0 \equiv 1\\
(2+1,-) & 2-0 \equiv 2\\
(1+1+1,-) & 3-0 \equiv 3\\
(2,1) & 1-1 \equiv 0\\
(1+1,1) & 2-1 \equiv 1\\
(1,2) & 1-1 \equiv 0\\
(1,1+1) & 1-2 \equiv4\\
(-,3) & 0-1 \equiv 4\\
(-,2+1) & 0-2 \equiv 3\\
(-,1+1+1) & 0-3 \equiv 2\\
\end{array}
\eeqs
Thus
$$
N_{\HL}(0,5,3) =  N_{\HL}(1,5,3) =  N_{\HL}(2,5,3) =  N_{\HL}(3,5,3) =  
N_{\HL}(4,5,3) =  2,
$$
and we see that the residue of the HL-birank mod $5$ divides the $10$
bipartitions of $3$ into $5$ equal classes.

\section{The Dyson-Birank} \label{sec:Dbirank}
Dyson \cite{Dy44}, \cite[p.52]{Dy96} defined the rank of a partition as the largest
part minus the number of parts. 
We define the Dyson-analog of the
birank for bipartitions $\bpi=\tpi$ by
\begin{equation}
\mbox{Dyson-birank}(\bpi) = \mbox{rank}(\pi_1)
+ 2\,\mbox{rank}(\pi_2).
\mylabel{eq:Dysonbirankdef2}
\end{equation}
In this section we prove
\begin{theorem}
\label{thm:thm2}
The residue of the Dyson-birank mod $5$
divides the bipartitions of $n$ into $5$ equal classes
provided $n\equiv2$, or $4\pmod{5}$.
\end{theorem}

We let $N_{\DY}(m,t,n)$ denote the number of bipartitions $\bpi=\tpi$
with Dyson-birank congruent to $m\pmod{t}$. 
We illustrate Theorem \thm{thm2} for the case $n=2$.
\beqs
\begin{array}{rr}
\mbox{Bipartitions of $2$} & \mbox{Dyson-birank (mod $5$)}\\
(2,-) & 1+0 \equiv 1\\
(1+1,-) & -1+0 \equiv 4\\
(1,1) & 0+0 \equiv 0\\
(-,2) & 0+2 \equiv 2\\
(-,1+1) & 0-2 \equiv 3\\
\end{array}
\eeqs
Thus
$$
N_{\DY}(0,5,2) =  N_{\DY}(1,5,2) =  N_{\DY}(2,5,2) =  N_{\DY}(3,5,2) =
N_{\DY}(4,5,2) =  1,
$$
and we see that the residue of the Dyson-birank mod $5$ divides the $5$
bipartitions of $2$ into $5$ equal classes. We note that Theorem
\thm{thm2} does not hold for $n\equiv 3\pmod{5}$. The first counterexample
occurs when $n=13$. The Dyson-birank mod $5$ fails to divide the
1770 bipartitions of $13$ into $5$ equal classes. We have
$$
N_{\DY}(0,5,13) =  358,\quad\mbox{but}\quad 
N_{\DY}(1,5,13) =  N_{\DY}(2,5,13) =  N_{\DY}(3,5,13) =
N_{\DY}(4,5,13) =  353.
$$

To prove Theorem \thm{thm2} we need the $5$-dissection of the rank generating
function when $z=\zeta$.
The Dyson-rank generating function is
\beq
f(z,q) = \sum_{\pi} z^{\srank(\pi)} q^{\abs{\pi}} =  1 +
\sum_{n=1}^\infty \frac{q^{n^2}}{(zq;q)_n (z^{-1}q;q)_n}.
\mylabel{eq:rankgenfunc}
\eeq
We let $N(m,t,n)$ denote the number of ordinary partitions of $n$
with rank congruent to $m$ mod $t$. Then
\begin{align}
f(\zeta,q)=&\sum_{n=0}^\infty 
\sum_{k=0}^4 \zeta^k N(k,5,n) q^n
=1+\sum_{n=1}^\infty \frac{q^{n^2}}{(\zeta q;q)_n ({\zeta}^{-1}q;q)_n}
\mylabel{eq:fdissect}\\
&= (A(q^5) - (3 + \zeta^2 + \zeta^3)\,\phi(q^5))
   + q\,B(q^5)
   + q^2 (\zeta + \zeta^4)\, C(q^5) \nonumber\\
&\qquad      
    +q^3 ( (1+\zeta^2+\zeta^3)\, D(q^5) + (1 + 2\zeta^2 + 2\zeta^3)\,\psi(q^5)),
\nonumber
\end{align}
where
\begin{align}
A(q) &=  \frac{E^2(q)\,J_{2,5}(q)}{J_{1,5}^2(q)},\mylabel{eq:Adef}\\
B(q) &=  \frac{E^2(q)}{J_{1,5}(q)},\mylabel{eq:Bdef}\\
C(q) &=  \frac{E^2(q)}{J_{2,5}(q)},\mylabel{eq:Cdef}\\
D(q) &=  \frac{E^2(q)\,J_{1,5}(q)}{J_{2,5}^2(q)},\mylabel{eq:Ddef}\\
\phi(q) &=-1+\sum_{n=0}^\infty\frac{q^{5n^2}}{(q;q^5)_{n+1}(q^4;q^5)_n} 
\mylabel{eq:phidef}\\
        &= \frac{q}{E(q^5)}\sum_{m=-\infty}^\infty
        (-1)^m\frac{q^{\frac{15}2m(m+1)}}{1-q^{5m+1}},\nonumber \\
\noalign{\mbox{and}}
\psi(q)
 &=\dfrac1{q}\biggl\{-1+\sum_{n=0}^\infty\frac{q^{5n^2}}{(q^2;q^5)_{n+1}(q^3;q^5)_n}
\biggr\} \mylabel{eq:psidef}\\
        &=\frac{q}{E(q^5)} \sum_{m=-\infty}^\infty
        (-1)^m\frac{q^{\frac{15}2m(m+1)}}{1-q^{5m+2}}.
\nonumber
\end{align}
Equation \eqn{fdissect} has an unusual history. It is originally due
to Ramanujan since appears in the Lost Notebook. It is closely related
to Dyson's conjectures on the rank \cite{Dy44}, which were proved
by Atkin and Swinnerton-Dyer \cite{At-Sw}. As pointed out in \cite{Ga88} and 
\cite{Ga88b},
equation \eqn{fdissect} is actually equivalent to one of
Atkin and Swinnerton-Dyer's main results. Dyson, Atkin and Swinnerton-Dyer
were unaware of Ramanujan's result.

The Dyson-birank generating function is
\beq
\sum_{\bpi=(\pi_1,\pi_2)} z^{\DYbr} q^{\abs{\bpi}} = f(z,q) \, f(z^2,q),
\mylabel{eq:DYgenfunc}
\eeq
where $f(z,q)$ is the generating function for the Dyson rank of ordinary partitions
given in \eqn{rankgenfunc}.
Thus we have
\beq
\sum_{n=0}^\infty \sum_{k=0}^4 \zeta^k N_{\DY}(k,5,n) q^n
= f(\zeta,q)\, f(\zeta^2,q).
\mylabel{eq:Dbrmod5}
\eeq
Using only \eqn{fdissect} and the fact that
\beq
B^2(q) = A(q)\,C(q),\qquad C^2(q) = B(q)\,D(q),
\mylabel{eq:ABCDfact}
\eeq
we find that the coefficient of $q^n$ in the $q$-expansion of
$f(\zeta,q)\, f(\zeta^2,q)$ is zero if 
$n\equiv2$, or $4\pmod{5}$. Theorem \thm{thm2} then follows from \eqn{Dbrmod5}.
Although Theorem \thm{thm2} does not hold when $n\equiv 3\pmod{5}$,
there is some simplification in the product $f(\zeta,q)\, f(\zeta^2,q)$.
We find that
\beq
\sum_{n=0}^\infty \sum_{k=0}^4 \zeta^k N_{\DY}(k,5,5n+3) q^n
= 5 \, \phi(q)\, \psi(q),         
\mylabel{eq:phipsiid}
\eeq
using the fact that
\beq
A(q) \, D(q) = B(q) \, C(q).
\mylabel{eq:ABCDfact2}
\eeq

\section{The $5$-Core-Birank} \label{sec:5cbirank}

For ordinary partitions $\pi$ the {\it $5$-core-crank} is defined
by
\begin{equation}
\mbox{$5$-core-crank}(\pi) = r_1 + 2 r_2 - 2 r_3 - r_4,
\mylabel{eq:GKScrankdef2}
\end{equation}
where $r_j$ is the number of cells labelled $j$ in the $5$-residue
diagram of $\pi$. See \cite[Prop.1,p.7]{Ga-Ki-St}. We define the 
$5$-core-crank
analog for bipartitions $\bpi=\tpi$ by
\begin{equation}
\mbox{$5$-core-birank}(\bpi) = \mbox{$5$-core-crank}(\pi_1)
+ 2\,(\mbox{$5$-core-crank}(\pi_2)).
\mylabel{eq:GKSbirankdef2}
\end{equation}

In this section we prove
\begin{theorem}
\label{thm:thm3}
The residue of the $5$-core-birank mod $5$
divides the bipartitions of $n$ into $5$ equal classes
provided $n\equiv2$, $3$ or $4\pmod{5}$.
\end{theorem}

We let $N_{\FC}(m,t,n)$ denote the number of bipartitions $\bpi=\tpi$
with $5$-core-birank congruent to $m\pmod{t}$. 
We illustrate Theorem \thm{thm3} for the case $n=3$.
\beqs
\begin{array}{rr}
\mbox{Bipartitions of $3$} & \mbox{$5$-core-birank (mod $5$)}\\
(3,-) & 3+0 \equiv 3\\
(2+1,-) & 0+0 \equiv 0\\
(1+1+1,-) & -3+0 \equiv 2\\
(2,1) & 1+0 \equiv 1\\
(1+1,1) & -1+0 \equiv 4\\
(1,2) & 0+2 \equiv 2\\
(1,1+1) & 0-2 \equiv3\\
(-,3) & 0+6 \equiv 1\\
(-,2+1) & 0+0 \equiv 0\\
(-,-6) & 0-3 \equiv 4\\
\end{array}
\eeqs
Thus
$$
N_{\FC}(0,5,3) =  N_{\FC}(1,5,3) =  N_{\FC}(2,5,3) =  N_{\FC}(3,5,3) =
N_{\FC}(4,5,3) =  2,
$$
and we see that the residue of the $5$-core-birank mod $5$ divides the $10$
bipartitions of $3$ into $5$ equal classes. We note that 
although the Dyson-birank does not in general divide the
bipartitions of $5n+3$ into $5$ equal classes the $5$-core-birank does.

To prove Theorem \thm{thm3} we need the $5$-dissection of the 
$5$-core-crank generating
function when $z=\zeta$.
The $5$-core-crank generating function is
\beq
\Phi(z,q) = \sum_{\pi} z^{\fccrank(\pi)} q^{\abs{\pi}} =  \frac{1}{E^5(q^5)}
T(z,q),
\mylabel{eq:fccrankgenfunc}
\eeq
where
\beq                  
T(z,q) := \sum_{\mbox{$\pi$ a $5$-core}}
z^{\mbox{$5$-core-crank($\pi$)}} q^{\abs{\pi}}
= \sum_{\substack{\vec{n}\in\Z^5\\\vec{n}\cdot\vec{1}=0}}
  z^{n_1 + 3n_2 + n_3} q^{\tfrac{5}{2}||\vec{n}||^2 + \vec{b}\cdot\vec{n}},
\mylabel{eq:Tdef}
\eeq
where $\vec{1}=(1,1,1,1,1)$ and $\vec{b}=(0,1,2,3,4)$.
Equation \eqn{fccrankgenfunc} can be proved combinatorially and
in a straightforward manner using Bijections 1 and 2 from 
\cite[pp.2-3]{Ga-Ki-St} and \cite[(4.2),  p.6]{Ga-Ki-St}.

We need the $5$-dissection of $T(\zeta,q)$:
\begin{equation}
T(\zeta,q) =  W(q^5)(1 + q R(q^5) + q^2 (\zeta^2+\zeta^3) R(q^5)^2
                          - q^3 (\zeta^2+\zeta^3) R(q^5)^3),
\mylabel{eq:T5dissect}
\end{equation}
where
\begin{align}
W(q) &:= 
J_{2,5}(q)^3( J_{10,25}(q) - q (1 + \zeta^2 + \zeta^3) J_{5,25}(q) ),
\mylabel{eq:Wdef}\\
\noalign{\mbox{and}}
R(q) &:= \frac{ J_{1,5}(q)}{J_{2,5}(q)}.\mylabel{eq:Rdef}
\end{align}
We will prove \eqn{T5dissect} in the next section. Theorem \thm{thm3}
follows easily from \eqn{T5dissect}.

The $5$-core-birank generating function is
\beq
\sum_{\bpi=(\pi_1,\pi_2)} z^{\fcbr} q^{\abs{\bpi}} = 
\frac{1}{E^{10}(q^5)}\, T(z,q) \, T(z^2,q),
\mylabel{eq:fcbrgenfunc}
\eeq
where $T(z,q)$ is the generating function for the 
$5$-core-crank of partitions that are $5$-cores given in \eqn{Tdef}.
Thus we have
\beq
\sum_{n=0}^\infty \sum_{k=0}^4 \zeta^k N_{\FC}(k,5,n) q^n
= \frac{1}{E^{10}(q^5)}\,T(\zeta,q)\, T(\zeta^2,q).
\mylabel{eq:FCbrmod5}
\eeq
From \eqn{T5dissect} we find that
\beq
T(\zeta,q) \, T(\zeta^2,q) = W^2(q^5)( 1 + 2\,q^5\,R^5(q^5)
+ q R(q^5)( 2 - q^5\,R^5(q^5)).
\mylabel{eq:TTid}
\eeq
Since coefficient of $q^n$ in the $q$-expansion of 
$T(\zeta,q) \, T(\zeta^2,q)$ is zero when 
$n\equiv2$, $3$ or $4\pmod{5}$, Theorem \thm{thm3} then follows
from \eqn{FCbrmod5}.

\section{A Theta-Function Identity} \label{sec:thetaid}

In this section we will prove the following theta-function identity.
\beq
U(z,q) = F_0(q)\, S_0(z,q) + F_1(q)\, S_1(z,q) + F_2(q)\, S_2(z,q)
+ F_3(q)\, S_3(z,q) + F_4(q)\, S_4(z,q),
\mylabel{eq:Uid}
\eeq
where
\begin{align}
F_0(q) &=  W(q^{10})\,(1 + q^{2}\,R(q^{10}) + (\zeta^2+\zeta^3)\,q^{4}\,R(q^{10})^2 - 
(\zeta^2+\zeta^3)\,q^{6}\,R(q^{10})^3),\mylabel{eq:f0def}\\
F_1(q) &=  W(q^{10})\,(\zeta^4 + \zeta\,q^{2}\,R(q^{10}) + (1+\zeta)\,q^{4},R(q^{10})^2 - 
(\zeta^2+\zeta^3)\,q^{6}\,R(q^{10})^3),\mylabel{eq:f1def}\\
F_2(q) &=  W(q^{10})\,(1 + \zeta^4\,q^{2}\,R(q^{10}) + (1+\zeta)\,q^{4}\,R(q^{10})^2 - 
(1+\zeta^4)\,q^{6}\,R(q^{10})^3),\mylabel{eq:f2def}\\
F_3(q) &=  W(q^{10})\,(1 + \zeta\,q^{2}\,R(q^{10}) + (1+\zeta^4)\,q^{4}\,R(q^{10})^2 - 
(1+\zeta)\,q^{6}\,R(q^{10})^3),\mylabel{eq:f3def}\\
F_4(q) &=  W(q^{10})\,(\zeta + \zeta^4\,q^{2}\,R(q^{10}) + (1+\zeta^4)\,q^{4}\,R(q^{10})^2 - 
(\zeta^2+\zeta^3)\,q^{6}\,R(q^{10})^3),\mylabel{eq:f4def}
\end{align}
\begin{align}
S_0(z,q) &= \sum_{n=-\infty}^{\infty} z^{5n}q^{25n^2+20n},
\mylabel{eq:S0def} \\
S_1(z,q) &= \sum_{n=-\infty}^{\infty} z^{5n+1}q^{25n^2+30n+5},
\mylabel{eq:S1def} \\
S_2(z,q) &= \sum_{n=-\infty}^{\infty} z^{5n+2}q^{25n^2+40n+12},
\mylabel{eq:S2def} \\
S_3(z,q) &= \sum_{n=-\infty}^{\infty} z^{5n+3}q^{25n^2+50n+21},
\mylabel{eq:S3def} \\
S_4(z,q) &= \sum_{n=-\infty}^{\infty} z^{5n+4}q^{25n^2+60n+32},
\mylabel{eq:S4def} 
\end{align}
and
\beq
U(z,q) 
 = \sum_{\substack{\vec{n}\in\Z^5}}
   z^{\vec{n}\cdot\vec{1}}
   \zeta^{n_1 + 3n_2 + n_3} 
   q^{5||\vec{n}||^2 + 2\vec{b}\cdot\vec{n}},
\mylabel{eq:Udef}
\eeq
where $W(q)$ and $R(q)$ defined in \eqn{Wdef} and \eqn{Rdef} respectively,
and the vectors 
$\vec{n}=(n_0,n_1,n_2,n_3,n_4$, $\vec{1}=(1,1,1,1,1)$ 
and $\vec{b}=(0,1,2,3,4)$ as before.
We note that \eqn{T5dissect} follows from \eqn{Uid} by taking the
coefficient of $z^0$ and replacing $q$ by $q^{1/2}$. Equation
\eqn{T5dissect} was the crucial identity needed in the proof of
Theorem \thm{thm3}.

We prove the identity \eqn{Uid} using standard techniques.
We show that both sides satisfy the same functional equation and
both sides agree for enough values of the parameter $z$. 
Most of these evaluations can be proved by elementary means using
Jacobi's triple product. For one evaluation we will need the theory
of modular functions.


We define the following Jacobi theta function
\beq
\Theta(z,q) = \sum_{n=-\infty}^\infty z^n q^{n^2},
\mylabel{eq:Thetadef}
\eeq
for $z\ne0$ and $\abs{q}<1$. We will need Jacobi's triple product
identity
\beq
\sum_{n=-\infty}^\infty z^n q^{n^2} 
=(-zq;q^2)_\infty (-z^{-1}q;q^2)_\infty (q^2;q^2)_\infty,
\mylabel{eq:jtprod}
\eeq
and the well-known functional equation
\beq
\Theta(zq^2,q) = z^{-1} q^{-1} \Theta(z,q),
\mylabel{eq:Thetafid}
\eeq
for $z\ne0$ and $0<\abs{q}<1$.
From the definition \eqn{Udef} we have
\beq
U(z,q) = \Theta(z\zeta^4,q^5) \, \Theta(zq^2\zeta,q^5) \, \Theta(zq^4,q^5)
        \, \Theta(z\zeta q^6,q^5) \, \Theta(z\zeta^4 q^8,q^5).
\mylabel{eq:Uprodid}
\eeq
From \eqn{jtprod} and \eqn{Thetafid} we have
\beq
U(z q^{10},q) = z^{-5} q^{-45} U(z,q),
\mylabel{eq:Ufid}
\eeq
and
\beq
U(z,q) = 0\quad\mbox{for}\quad  
z = -q^5\,\zeta,\quad  -q^3\,\zeta^4, -q, -\zeta^4\,q^9, -\zeta\,q^7.
\mylabel{eq:Uzeros}
\eeq
Let $V(z,q)$ denote the function of the right side of \eqn{Uid}.
Each $S_j(z,q)$ can be written in terms of the theta function
$\Theta(z,q)$ and we find that 
$S_j(z q^{10},q) = z^{-5} q^{-45} S_j(z,q)$ for each $j$ so that
\beq
V(z q^{10},q) = z^{-5} q^{-45} V(z,q).
\mylabel{eq:Vfid}
\eeq
Hence the left and right sides of \eqn{Uid} satisfy the same
functional equation (i.e.\ \eqn{Ufid}, \eqn{Vfid}).
In view of \cite[Lemma 2]{At-Sw} or \cite[Lemma 1]{Hi-Ga-Bo},
it suffices to show that \eqn{Uid} holds for $6$ distinct
values of $z$ with $\abs{q}^{10} < \abs{z} \le 1$.
We claim that
\beq
V(z,q) = 0\quad\mbox{for}\quad  
z = -q^5\,\zeta,\quad  -q^3\,\zeta^4, -q, -\zeta^4\,q^9, -\zeta\,q^7.
\mylabel{eq:Vzeros}
\eeq
Using \eqn{jtprod} we can easily evaluate each $S_j(z,q)$ for
these values of $z$.
\begin{align}
  S_0(-\zeta q^5,q) &= -q^{-20}J_{2,5}(q^{10}),
\mylabel{eq:Sv1}\\
  S_1(-\zeta q^5,q) &= \zeta q^{-20} J_{2,5}(q^{10}),
\mylabel{eq:Sv2}\\
  S_2(-\zeta q^5,q) &= -\zeta^2 q^{-18} J_{1,5}(q^{10}),
\mylabel{eq:Sv3}\\
  S_3(-\zeta q^5,q) &= 0,
\mylabel{eq:Sv4}\\
  S_4(-\zeta q^5,q) &= \zeta^4 q^{-18} J_{1,5}(q^{10}),
\mylabel{eq:Sv5}\\
\noalign{\medskip}
  S_0(-q^3\zeta^4,q) &= -q^{-10}J_{1,5}(q^{10}),
\mylabel{eq:Sv6}\\
  S_1(-q^3\zeta^4,q) &= \zeta^4 q^{-12} J_{2,5}(q^{10}),
\mylabel{eq:Sv7}\\
  S_2(-q^3\zeta^4,q) &= -\zeta^3 q^{-12} J_{2,5}(q^{10}),
\mylabel{eq:Sv8}\\
  S_3(-q^3\zeta^4,q) &= \zeta^2 q^{-10} J_{1,5}(q^{10}),
\mylabel{eq:Sv9}\\
  S_4(-q^3\zeta^4,q) &= 0 ,                                       
\mylabel{eq:Sv10}\\
\noalign{\medskip}
  S_0(-q,q) &= 0,
\mylabel{eq:Sv11}\\
  S_1(-q,q) &= q^{-4}J_{1,5}(q^{10}),
\mylabel{eq:Sv12}\\
  S_2(-q,q) &= -q^{-6}J_{2,5}(q^{10}),
\mylabel{eq:Sv13}\\
  S_3(-q,q) &= q^{-6}J_{2,5}(q^{10}),
\mylabel{eq:Sv14}\\
  S_4(-q,q) &= -q^{-4}J_{1,5}(q^{10}),
\mylabel{eq:Sv15}\\
\noalign{\medskip}
  S_0(-\zeta^4q^9,q) &= -q^{-40}J_{1,5}(q^{10}),
\mylabel{eq:Sv16}\\
  S_1(-\zeta^4q^9,q) &= 0
\mylabel{eq:Sv17}\\
  S_2(-\zeta^4q^9,q) &= \zeta^3q^{-40}J_{1,5}(q^{10}),
\mylabel{eq:Sv18}\\
  S_3(-\zeta^4q^9,q) &= -\zeta^2q^{-42}J_{2,5}(q^{10}),
\mylabel{eq:Sv19}\\
  S_4(-\zeta^4q^9,q) &= \zeta q^{-42}J_{2,5}(q^{10}),
\mylabel{eq:Sv20}\\
\noalign{\medskip}
  S_0(-\zeta q^7,q) &= -q^{-30}J_{2,5}(q^{10}),
\mylabel{eq:Sv21}\\
  S_1(-\zeta q^7,q) &= \zeta q^{-28}J_{1,5}(q^{10}),
\mylabel{eq:Sv22}\\
  S_2(-\zeta q^7,q) &= 0,
\mylabel{eq:Sv23}\\
  S_3(-\zeta q^7,q) &= -\zeta^3q^{-28}J_{1,5}(q^{10}),
\mylabel{eq:Sv24}\\
  S_4(-\zeta q^7,q) &= \zeta^4q^{-30}J_{2,5}(q^{10}).
\mylabel{eq:Sv25}
\end{align}
The verification of \eqn{Vzeros} is just a routine calculation.

Thus both sides of \eqn{Uid} agree for $5$ distinct values of $z$
in the region $\abs{q}^{10} < \abs{z} \le 1$.
We show that both sides agree for $z=-1$, and then our
identity \eqn{Uid} will follow. To achieve this we use the theory
of modular functions. Since this is a standard technique we just
sketch some of the details.

First, we calculate the $5$-dissection of each theta function on 
the right side of \eqn{Uprodid} when $z=-1$. By \eqn{jtprod} we find
that
\begin{align}
\Theta(-\zeta^4,q^5) &= J_{1,2}(q^{50})+(1+\zeta^2+\zeta^3)q^{5}J_{3,10}(q^{25})
            +(\zeta^2+\zeta^3)q^{20}J_{1,10}(q^{25}),
\mylabel{eq:Peval1}\\
\Theta(-\zeta q^2,q^5)&= J_{27,50}(q^{5})+\zeta^3q^{16}J_{7,50}(q^{5})
            -\zeta q^{7}J_{13,50}(q^{5})-\zeta^4q^{3}J_{17,50}(q^{5})
            +q^{24}\zeta^2J_{7,50}(q^{5}),
\mylabel{eq:Peval2}\\
\Theta(-q^4,q^5) &= J_{43,50}(q^{5})-q\,J_{19,50}(q^{5})+q^{12}J_{9,50}(q^{5})
            +q^{28}J_{1,50}(q^{5}) -q^{9}J_{11,50}(q^{5}),
\mylabel{eq:Peval3}\\
q\,\Theta(-\zeta q^6,q^5)&= -\zeta^4J_{21,50}(q^{5})+q\,J_{19,50}(q^{5})
               -\zeta q^{12}J_{9,50}(q^{5})
               -\zeta^2q^{28}J_{1,50}(q^{5}) +\zeta^3q^{9}J_{23,50}(q^{5}),
\mylabel{eq:Peval4}\\
q^{3} \Theta(-\zeta^4 q^8,q^5)
&=-\zeta J_{23,50}(q^{5})-\zeta^4q^{16}J_{7,50}(q^{5})
               +\zeta^2q^{7}J_{27,50}(q^{5}) +q^{3}J_{17,50}(q^{5})
               -\zeta^3q^{24}J_{3,50}(q^{5}).
\mylabel{eq:Peval5}
\end{align}
Next, we evaluate each $S_j(-1,q)$ using \eqn{jtprod}
\begin{align}
 S_0(-1,q) &=  S_1(-1,q) = J_{1,10}(q^{5}),\mylabel{eq:Seval1}\\
 S_2(-1,q) &=  S_4(-1,q) = -q^{-3}J_{3,10}(q^{5}),\mylabel{eq:Seval2}\\
 S_3(-1,q) &= q^{-4}J_{1,2}(q^{25}).\mylabel{eq:Seval3}          
\end{align}

For $0\le r \le 4$, we define the operator $\mathcal{U}_{r,5}$ by
\beq
\mathcal{U}_{r,5}\left(\sum_{n} a(n) q^n \right)
= \sum_{n} a(5n+r) q^n.
\mylabel{eq:Uop}
\eeq
To show that \eqn{Uid} holds for $z=-1$ we need to prove $5$ identities
\beq
\mathcal{U}_{r,5}\left(q^4\,U(-1,q)\right)
=
\mathcal{U}_{r,5}\left(q^4\,V(-1,q)\right),
\mylabel{eq:UVrid}
\eeq
for $0\le r \le 4$. It turns out that each of these identities is equivalent
to a modular function identity for the group $\Gamma_1(50)$. We provide
some detail for the case $r=0$. Using \eqn{Seval1}--\eqn{Seval3} we 
find that 
\begin{align}
\mathcal{U}_{0,5}\left(q^4\,V(-1,q)\right) 
&=
\left(J_{2,5}(q^{10}) - q^2 (1 + \zeta^2 + \zeta^3) J_{1,5}(q^{10})\right)
\mylabel{eq:U0Vid}\\
&\times \left( J_{1,2}(q^5) J_{4,10}^3(q)
       + q (-1 + \zeta^2 + \zeta^3) J_{2,10}^2(q) J_{3,10}(q) J_{4,10}(q)\right. \nonumber\\
&\hphantom{XXXXXXXXXXXXXX}\left. - 2 q^2 (\zeta^2+\zeta^3) J_{2,10}^3 J_{1,10}(q) \right).
\nonumber
\end{align}
We can utilize \eqn{Peval1}--\eqn{Peval5} to write the left side
of the $r=0$ case of \eqn{UVrid} as a sum of $135$ explicit theta products
\begin{align}
&\mathcal{U}_{0,5}\left(q^4\,U(-1,q)\right) \mylabel{eq:U0Uid}\\
&\quad =
q^{10} J_{3,50}^2(q) J_{19,50}^2(q) J_{25,50}(q) + \cdots
+ 2 q^{10} (\zeta^2 + \zeta^3) J_{1,50}(q) J_{9,50}(q) J_{13,50}(q) J_{17,50}(q)
  J_{25,50}(q).
\nonumber
\end{align}
We have to prove that the right side of \eqn{U0Vid} equals the right side of \eqn{U0Vid}.
After dividing both sides by $J_{2,5}(q^{10}) J_{1,2}(q^5) J_{4,10}^3(q)$ we find that this
is equivalent to showing that a certain linear combination of  $140$ generalized 
eta-quotients 
simplifies to the constant $1$   
\beq
(1+\zeta^2+\zeta^3)\,\eta_{50,10} \, \eta_{50,20}^{-1}
+ \cdots
\eta_{50,4}^{-3} \, \eta_{50,5}^{-2} \, \eta_{50,6}^{-3} \, \eta_{50,10}^{-4} \,
\eta_{50,14}^{-3} \, \eta_{50,15}^{-2} \, \eta_{50,16}^{-3} \, \eta_{50,20}^{-5} \,
\eta_{50,24}^{-3} \, \eta_{50,23}^{2} \, \eta_{50,21}^{2} = 1.
\mylabel{eq:bigcombo}
\eeq
Here
\beq
\eta_{n,m} = \eta_{n,m}(\tau)
= \exp(\pi i P_2(m/n) n\tau)\,
  \prod_{k\equiv\pm m\pmod{n}}(1 - \exp(2\pi ik\tau)
= q^{n P_2(m/n)/2} J_{m,n}(q),
\mylabel{eq:getadef}
\eeq
where $P_2(t) = \{t\}^2 - \{t\} + \tfrac{1}{6}$, and $q=\exp(2\pi i\tau)$.
Using \cite[Theorem 2.9, p.7]{Choi06}, \cite[Theorem 3, p.126]{Ro94} that
each generalised eta-quotient in \eqn{bigcombo} is indeed a modular
function on $\Gamma_1(50)$. As usual we need the  valence formula
\beq
\sum_{z\in\mathcal{F}} \ORD(f;z,\Gamma) = 0,
\mylabel{eq:valform}
\eeq
provided $f$ is a nontrivial modular function on $\Gamma$, and
$\mathcal{F}$ is a fundamental set for $\Gamma$. Using MAGMA, the following
is a complete set of inequivalent cusps for $\Gamma_1(50)$
\begin{align}
\mathcal{C}
&=\{\infty, \tfrac{0}{1},\, \tfrac{1}{10},\, \tfrac{1}{9},\, \tfrac{2}{17},\, \tfrac{3}{25},\,
\tfrac{1}{8},\, \tfrac{2}{15},\, \tfrac{3}{22},\, \tfrac{4}{29},\, \tfrac{5}{36},\, 
\tfrac{6}{43},\, \tfrac{7}{50},\, \tfrac{3}{20},\, \tfrac{2}{13},\, \tfrac{3}{19},\, 
\tfrac{4}{25},\, \tfrac{1}{6},\, \tfrac{6}{35},\, \tfrac{9}{52},\, \tfrac{13}{75},\, 
\tfrac{4}{23},\, \tfrac{7}{39},\, \tfrac{9}{50},\, \tfrac{7}{38},\, \tfrac{3}{16},\, 
\mylabel{eq:CG}\\
&\qquad\tfrac{5}{26},\, \tfrac{1}{5},\, \tfrac{27}{125},\, \tfrac{7}{32},\, \tfrac{11}{50},\,
\tfrac{12}{53},\, \tfrac{17}{75},\, \tfrac{6}{25},\, \tfrac{1}{4},\, \tfrac{13}{50},\, 
\tfrac{53}{200},\, \tfrac{67}{250},\, \tfrac{27}{100},\, \tfrac{11}{40},\, \tfrac{18}{65},\, 
\tfrac{7}{25},\, \tfrac{3}{10},\, \tfrac{7}{20},\, \tfrac{9}{25},\, \tfrac{11}{30},\, 
\tfrac{19}{50},\, \tfrac{49}{125},\, \tfrac{2}{5},\, \tfrac{21}{50},\, \nonumber\\
&\qquad \tfrac{11}{25},\, 
\tfrac{11}{20},\, \tfrac{26}{45},\, \tfrac{3}{5},\, \tfrac{59}{85},\, \tfrac{7}{10}\},
\nonumber
\end{align}
with corresponding widths
\begin{align}
&\{
1,\, 50,\, 5,\, 50,\, 50,\, 2,\, 25,\, 10,\, 25,\, 50,\, 25,\, 50,\, 1,\, 5,\, 50,\, 50,\, 2,\, 25,\, 10,\, 25,\, 2,\, 50,\, 50,\, 1,\, 25,\,
\mylabel{eq:WG}\\
&\qquad    25,\, 25,\, 10,\, 2,\, 25,\, 1,\, 50,\, 2,\, 2,\, 25,\, 1,\, 1,\, 1,\, 1,\, 5,\, 10,\, 2,\, 5,\, 5,\, 2,\, 5,\, 1,\, 2,\, 10,\, 1,\, 2,\, 5,\, 
\nonumber\\
&\qquad   10,\, 10,\, 10,\, 5\}.
\nonumber
\end{align}
Using known results for the invariant order of generalized eta-quotients
at cusps \cite[(2.3), p.7]{Choi06}, \cite[pp.127-128]{Ro94} we have calculated the
order at each cusp of every function in \eqn{bigcombo}. As check we verified
that the total Order of each function is zero. With $\mathcal{J}$ being
the set generalized eta-quotients ocurring in \eqn{bigcombo} we calculated
\beq
\sum_{c\in\mathcal{C}\setminus\{\infty\}}
\min_{ f\in \mathcal{J}}(\ORD(f;c;\Gamma_1(50)),0) = -145.
\mylabel{eq:minords}
\eeq
Hence, by the valence formula \eqn{valform} it suffices to verify 
\eqn{bigcombo} (or equivalently \eqn{UVrid} with $r=0$) 
up to $q^{145}$, since generalized eta-quotients have no poles or zeros in the upper-half
plane. We have actually verified
the result up to $q^{200}$. All calculations, except for \eqn{CG},and \eqn{WG}, 
were done using MAPLE. The calculations needed to verify  
\eqn{UVrid} for $r=1$, $2$, $3$, $4$ are similiar and have been carried out.
This conpletes our proof of \eqn{Uid}.

\section{The Andrews Bicrank and Extensions} \label{sec:bicrank}

For a partition $\pi$, let $\ell(\pi)$ denote the largest part of $\pi$, 
$\varpi(\pi)$ denote the number of ones in $\pi$,
and $\mu(\pi)$ denote the number of parts of $\pi$ larger than $\varpi(\pi)$.
The crank of $\pi$ is given by
\beq
\mbox{crank}(\pi) =
\begin{cases}
  \ell(\pi), &\mbox{if $\varpi(\pi)=0$}, \\
  \mu(\pi) - \varpi(\pi), &\mbox{if $\varpi(\pi)>0$}.
\end{cases}
\mylabel{eq:crankdef}
\eeq
The crank gives a combinatorial interpretation
of Ramanujan's partition congruences mod $5$, $7$ and $11$ and solves a
problem of Dyson \cite{Dy44}, \cite[p.52]{Dy96}. See \cite{An-Ga88}. 

In \cite{An08}, Andrews gave a combinatorial interpretation of the
congruence
\beq
p_{-2}(5n+3) \equiv 0 \pmod{5},
\mylabel{eq:p25n3}
\eeq
in terms of the crank. This result is a crank analog of the Dyson-birank
but is more complicated since it on involves positive and negative weights.
This complication is because of the nature of the generating function
for the crank. Let $M(m,n)$ denote the number of partitions of $n$ with
crank $m$. Then
\beq
\sum_{n\ge0} \sum_m M(m,n) z^m q^n = (1-z)q + 
\frac{(q;q)_\infty}{(zq;q)_\infty (z^{-1}q;q)_\infty}.
\mylabel{eq:Mgf}
\eeq
Define $M'(m,n)$ by
\beq
\sum_{n\ge0} \sum_m M'(m,n) z^m q^n = 
\frac{(q;q)_\infty}{(zq;q)_\infty (z^{-1}q;q)_\infty} 
= 1 + (z -1 + z^{-1})q + (z^2 + z^{-2})q^2 + \cdots.
\mylabel{eq:Mvgf}
\eeq
We need to interpret $M'(m,n)$ combinatorially. To this we need to
definition of partition.
To the set of partitions we need to add two additional partitions of $1$ which we denote
by $1_a$ and $1_b$. We call this new set $\mathcal{E}$, the set of extended partitions.
\beq
\mathcal{E} = \{(-), 1_a, 1_b, 1, 2, 1+1, 3, 2+1, 1+1+1, \cdots\}.
\mylabel{eq:eptns}
\eeq
We have $\abs{1_a}=\abs{1_b}=1$.
Here as usual $(-)$ is the empty partition of $0$. For these extended partitions
define a weight function $w(\pi)$ defined by
\beq
w(\pi) =
\begin{cases}
  -1,   &\mbox{if $\pi=1_b$}.\\
  \abs{\pi},   &\mbox{otherwise}.
\end{cases}
\mylabel{eq:wdef}
\eeq
Thus for the three extended partitions of $1$ we have
$w(1)=w(1_a)=1$, and $w(1_b)=-1$, and the total weight is still $p(1)=1$.
Therefore
\beq
\sum_{\substack{\pi\in\mathcal{E}\\ \abs{\pi}=n}}  w(\pi) = p(n).
\mylabel{eq:totw}
\eeq
We also extend the definition of crank by $\mbox{crank}(1_a)=1$,
and $\mbox{crank}(1_b)=0$. Recall that for ordinary partition of $1$ we have
$\mbox{crank}(1_b)=-1$.
We now have our desired combinatorial interpretation of $M'(m,n)$.
\beq
F(z,q) = \sum_{\pi\in\mathcal{E}} w(\pi) z^{\mbox{crank}(\pi)} q^{\abs{\pi}} 
= \sum_{n\ge0} \sum_m M'(m,n) z^m q^n = 
\frac{(q;q)_\infty}{(zq;q)_\infty (z^{-1}q;q)_\infty}.
\mylabel{eq:Mvid}
\eeq
In other words,
\beq
M'(m,n) = 
\sum_{\substack{\pi\in\mathcal{E}\\ \abs{\pi}=n,\, \mbox{crank}(\pi)=m}}  w(\pi).
\mylabel{eq:Mvid2}
\eeq
We note that the function $F(z,q)$ (at least as an infinite product)
occured in Ramanujan's Lost Notebook.

We define the set of extended bipartitions by $\mathcal{E} \times \mathcal{E}$,i.e.\ an
extended bipartition is simply a pair of extended partitions. For
an extended bipartition $\pi=(\pi_1,\pi_2)$ we define
a sum of parts function and a weight function in the natural way
\beq
\abs{\pi} = \abs{\pi_1} + \abs{\pi_2},\quad\mbox{and}\qquad
w(\pi) = w(\pi_1) \, w(\pi_2).
\mylabel{eq:ewdef}
\eeq
We denote Andrews's bicrank function by $\mbox{bicrank}_1$.
We give a variant which we call $\mbox{bicrank}_2$.
For an extended bipartition $\pi=(\pi_1,\pi_2)$ we define
\begin{align}
\mbox{bicrank}_1(\pi) &= \mbox{crank}(\pi_1) + \mbox{crank}(\pi_2),
\mylabel{eq:bic1}\\
\mbox{bicrank}_2(\pi) &= \mbox{crank}(\pi_1) + 2\,\mbox{crank}(\pi_2),
\mylabel{eq:bic2}
\end{align}
Amazingly together these two bicrank functions give a
new interpretation for all three congruences in \eqn{2colcongs}.
For $j=1$, $2$ we define  $M_j(m,t,n)$ by
\beq
M_j(m,t,n) =
\sum_{\substack{\pi\in\mathcal{E}\times\mathcal{E}\\ 
     \abs{\pi}=n,\, \mbox{bicrank}_j(\pi)\equiv m\pmod{t}}}  w(\pi).
\mylabel{eq:Mjdef}
\eeq                     
In other words, $M_j(m,t,n)$ is the number of extended bipartitions of $n$
with  $\mbox{bicrank}_j$ congruent to $m$ mod $t$ counted by the weight $w$.

In this section we prove
\begin{theorem}
\label{thm:thm4}
$$
\hphantom{x}
$$
\begin{enumerate}
\item[(i)]
The residue of the  $\mbox{bicrank}_1(\pi)$ mod $5$
divides the extended bipartitions of $n$ into $5$ classes of equal weight
provided $n\equiv3\pmod{5}$.
\item[(ii)]
The residue of the  $\mbox{bicrank}_2(\pi)$ mod $5$
divides the extended bipartitions of $n$ into $5$ classes of equal weight
provided $n\equiv2$ or $4\pmod{5}$.
\end{enumerate}
\end{theorem}

We illustrate the first case of Theorem \thm{thm4} (i).                         
There are $18$ extended bipartitions of $3$ giving a total weight of $p_{-2}(3)=10$.
\beqs
\begin{array}{rrr}
\mbox{Extended bipartitions of $3$} & \mbox{bicrank}_1 \pmod{5} & \mbox{weight}=w \\
(2+1,-)  & 0 \equiv 0    &  1\\
(-,2+1) &  0+0 \equiv 0  &  1\\
(2,1)    & 2-1 \equiv 1  &  1\\
(1,2)    &-1+2 \equiv 1  &  1\\
(1+1+1,-)& -3 \equiv 2   &  1\\
(2,1_b)  & 2+0 \equiv 2  & -1\\
(1+1,1)  & -2-1 \equiv 2 &  1\\
(1_b,2)  & 0+2 \equiv 2  & -1\\
(1,1+1)  &-1-2 \equiv 2   &  1\\
(-,1+1+1) &0-3 \equiv 2  &  1\\
(3,-)    & 3 \equiv 3    &  1\\
(2,1_a)  & 2+1 \equiv 3  &  1\\
(1+1,1_b)&-2+0 \equiv 3  & -1\\
(1_a,2)  & 1+2 \equiv 3  &  1\\
(1_b,1+1)& 0-2 \equiv 3   & -1\\
(-,3)    & 0+3 \equiv 3  &  1\\
(1+1,1_a)& -2+1 \equiv 4 &  1\\
(1_a,1+1)& 1-2 \equiv 4   &  1\\
\end{array}
\eeqs
Thus
$$
M_1(0,5,3) =  M_1(1,5,3) =  M_1(2,5,3) =  M_1(3,5,3) =  
M_1(4,5,3) =  2,
$$
and we see that the residue of the $\mbox{bicrank}_1$  mod $5$ divides the $18$
bipartitions of $3$ into $5$ classes of equal total weight $2$.

We illustrate the first case of Theorem \thm{thm4} (ii).                         
There are $13$ extended bipartitions of $2$ giving a total weight of $p_{-2}(2)=5$.
\beqs
\begin{array}{rrr}
\mbox{Extended bipartitions of $2$} & \mbox{bicrank}_2 \pmod{5} & \mbox{weight}=w \\
(2,-)     & 2+0 \equiv 2     & 1\\
(1+1,-)   &-2+0 \equiv 3     & 1\\
(1,1) & -1-2 \equiv 2        & 1\\
(1,1_a) & -1+2 \equiv 1      & 1\\
(1,1_b) & -1+0 \equiv 4      &-1\\
(1_a,1) & 1-2 \equiv 4       & 1\\
(1_a,1_a)& 1+2 \equiv 3      & 1\\
(1_a,1_b) & 1+0 \equiv 1     &-1\\
(1_b,1) & 0-2 \equiv 3       &-1\\
(1_b,1_a) & 0+2 \equiv 2     &-1\\
(1_b,1_b) & 0+0 \equiv 0     & 1\\
(-,2) & 0+4 \equiv 4         & 1\\
(-,1+1) & 0-4 \equiv 1       & 1\\
\end{array}
\eeqs
Thus
$$
M_2(0,5,2) =  M_2(1,5,2) =  M_2(2,5,2) =  M_2(3,5,2) =  
M_2(4,5,2) =  1,
$$
and we see that the residue of the $\mbox{bicrank}_2$  mod $5$ divides the $13$
bipartitions of $2$ into $5$ classes of equal total weight $1$.

Theorem \thm{thm4} (i) is due to Andrews \cite{An08}. 
Theorem \thm{thm4} (ii) is a natural extension, and its proof is analogous.
For completeness we include a sketch of the proof. We define
\beq
M'(r,t,n) = \sum_{m\equiv r\pmod{t}} M'(m,n),
\eeq
which is the number of ordinary partitions of $n$ with crank congruent
to $r$ mod $t$ when $n\ne1$. When $n=1$ it is counting extended partitions.
Then
\begin{align}
F(\zeta,q)=&\sum_{n=0}^\infty
\sum_{k=0}^4 \zeta^k M'(k,5,n) q^n
=  \frac{(q;q)_\infty}{(\zeta q;q)_\infty (\zeta^{-1}q;q)_\infty}
\mylabel{eq:Fdissect}\\
&= A(q^5)
   - q (\zeta + \zeta^4)^2 \,B(q^5)
   + q^2 (\zeta^2 + \zeta^3)\, C(q^5)  -q^3 (\zeta+\zeta^4)\, D(q^5),
\nonumber
\end{align}
where
$F(z,q)$ is given in \eqn{Mvid}, $A(q)$, $B(q)$, $C(q)$, and $D(q)$ are
given in \eqn{Adef}--\eqn{Ddef}. Equation \eqn{Fdissect} appears
in Ramanujan's Lost Notebook \cite[p.20]{Ra88} and is
proved in \cite[(1.30)]{Ga88}.

The two bicrank generating functions are given by
\begin{align}
\sum_{\bpi=(\pi_1,\pi_2)} z^{\bcra} q^{\abs{\bpi}} &= F(z,q)^2, 
\mylabel{eq:bcr1genfunc}\\
\sum_{\bpi=(\pi_1,\pi_2)} z^{\bcrb} q^{\abs{\bpi}} &= F(z,q) \, F(z^2,q),
\mylabel{eq:bcr2genfunc}
\end{align}
where $F(z,q)$ is the generating function for the crank of
extended partitions \eqn{Mvid}.
Thus we have
\begin{align}
\sum_{n=0}^\infty \sum_{k=0}^4 \zeta^k M_{1}(k,5,n) q^n
&= F(\zeta,q)^2,
\mylabel{eq:bcr1mod5}\\
\sum_{n=0}^\infty \sum_{k=0}^4 \zeta^k M_{1}(k,5,n) q^n
&= F(\zeta,q)\, F(\zeta^2,q).          
\mylabel{eq:bcr2mod5}
\end{align}
Using only \eqn{Fdissect} and equations
\eqn{ABCDfact} and \eqn{ABCDfact2}
we easily find that 
find that the coefficient of $q^n$ in the $q$-expansion of
$F(\zeta,q)^2$  is zero if
$n\equiv3\pmod{5}$, ant that the 
the coefficient of $q^n$ in the $q$-expansion of
$F(\zeta,q)\, F(\zeta^2,q)$ is zero if
$n\equiv2$, or $4\pmod{5}$. Both parts of Theorem \thm{thm4} then follow
from equations \eqn{bcr1mod5} and \eqn{bcr2mod5}.

\section{A Multirank Analog of the Hammond-Lewis Birank} \label{sec:ghlbirank}

Let $\mathcal{P}$ denote the set of partitions. A multipartition 
with $r$ components or an $r$-colored partition of $n$ is simply an
$r$-tuple
\beq
\vec{\pi}=(\pi_1,\pi_2,\dots,\pi_r) \in
\mathcal{P} \times \mathcal{P} \times \cdots \times \mathcal{P} 
=\mathcal{P}^r,
\eeq
where
\beq
\sum_{k=1}^r \abs{\pi_k}=n.
\eeq
It is clear that the number of $r$-colored partitions of $n$ is
$p_{-r}(n)$ where
\beq
\sum_{n\ge0} p_{-r}(n) q^n = \frac{1}{E(q)^r}.
\eeq
There are two elementary and well-known congruences.

\begin{theorem}
\label{thm:thm5}

Let $t>3$ be prime.

\begin{enumerate}
\item[(i)] If 
           $24n+1$ is a quadratic nonresidue mod $t$, then
 then
\beq
p_{1-t}(n) \equiv 0 \pmod{t}.
\eeq
\item[(ii)] If 
             $8n+1$ is not a quadratic residue mod $t$,
             then
\beq
p_{3-t}(n) \equiv 0 \pmod{t}.
\eeq
\end{enumerate}
\end{theorem}

These results follow easily from identities of Euler and Jacobi.
Theorem \thm{thm5} (i) follows from
\beq
\sum_{n\ge0} p_{1-t}(n) q^{24n+1} = \frac{qE(q^{24})}{E(q^{24})^t}
\equiv \frac{1}{E(q^{24t})} \sum_{n=-\infty}^\infty (-1)^n q^{(6n+1)^2}
\pmod{t}.
\eeq
Here we have used Euler's Pentagonal Number Theorem \cite[Thm 353]{Ha-Wr-BOOK}
\beq
E(q) = \sum_{n=-\infty}^\infty (-1)^n q^{n(3n+1)/2}.
\mylabel{eq:Epent}
\eeq
Theorem \thm{thm5} (ii) follows from
\beq
\sum_{n\ge0} p_{3-t}(n) q^{8n+1} = \frac{qE^3(q^{8})}{E(q^{8})^t}
\equiv \frac{1}{E(q^{8t})} \sum_{n\ge0}(-1)^n (2n+1) q^{(2n+1)^2}
\pmod{t},
\eeq
where we have used Jacobi's Identity \cite[Thm 237]{Ha-Wr-BOOK}
\beq
E(q)^3 = \sum_{n\ge0} (-1)^n (2n+1) q^{n(n+1)/2}.
\mylabel{eq:jacid}
\eeq
Theorem \thm{thm6} (ii) is Theorem 1 in \cite{An08}.

In this section we construct analogs of the Hammond-Lewis birank to
combinatorially explain the two congruences in Theorem \thm{thm5}.
Andrews's bicrank \cite{An08} (see also equation \eqn{bic1}) gave
a combinatorial interpretation of Theorem \thm{thm5} (ii) for the
case $t=5$, and $n\equiv 3\pmod{5}$. The Hammond-Lewis birank gave 
a combinatorial interpretation of Theorem \thm{thm5} (ii) for
the case $t=5$, and all relelvant $n$. 


For even $r$, we define
the generalized Hammond-Lewis multirank by
\begin{equation}
\mbox{gHL-multirank}(\vec{\pi}) 
= \sum_{k=1}^{r/2} k\,\left(\#(\pi_k) - \#(\pi_{r+1-k})\right),
\mylabel{eq:gHLmrankdef}
\end{equation}
for $\vec{\pi}=(\pi_1,\pi_2,\dots,\pi_{r})$ a multipartition with
$r$ components. The $r=2$ case corresponds to the Hammond-Lewis birank.

In this section we prove
\begin{theorem}
\label{thm:thm6}

Let $t>3$ be prime. 

\begin{enumerate}
\item[(i)] 
The residue of the  generalized-Hammond-Lewis-multirank mod $t$
divides the multipartitions of $n$ with $r=t-1$ components into $t$ equal classes
provided $24n+1$ is a quadratic nonresidue mod $t$.                   
\item[(ii)] 
The residue of the  generalized-Hammond-Lewis-multirank mod $t$
divides the multipartitions of $n$ with $r=t-3$ components into $t$ equal classes
provided   $8n+1$ is not a quadratic residue mod $t$.
\end{enumerate}
\end{theorem}

We illustrate  Theorem \thm{thm6} (ii) for $t=7$ and $n=2$.                         
\beqs
\begin{array}{rr}
\mbox{Multipartitions of $2$} & \mbox{generalized-HL-multirank}\\
\mbox{with $4$ components} & \mbox{(mod $7$)}\\
(-, -, -, 1+1)    & -2\equiv5\\
( -, -, -, 2)      & -1\equiv6\\
( -, -, 1, 1)     & -3\equiv4\\
(-, -, 1+1,-)    & -4\equiv3\\
( -, -, 2,-)      & -2\equiv5\\
( -, 1, -, 1)     &  1\equiv1\\
( -, 1, 1,-)     &  0\equiv0\\
(-, 1+1, -,-)    &  4\equiv4\\
( -, 2, -,-)      &  2\equiv2\\
( 1, -, -, 1)     &  0\equiv0\\
( 1, -, 1,-)     & -1\equiv6\\
( 1, 1, -,-)     &  3\equiv3\\
(1+1, -, -,-)    &  2\equiv2\\
( 2, -, -,-)      &  1\equiv1\\
\end{array}
\eeqs
We see that the residue of generalized-Hammond-Lewis-multirank mod $7$
divides the $14$ $4$-colored partitions of $2$ into $7$ equal classes.

Both parts of Theorem \thm{thm6} are easy to prove. For (i),
we need only Euler's pentagonal number theorem \eqn{Epent}.
We let $\zeta_t$ be a primitive $t$-th root
of unity. We have
\begin{align}
\sum_{\vec{\pi}\in\mathcal{P}^{t-1}}  \zeta_t^{\gHLmr} q^{\abs{\vec{\pi}}} 
&= \prod_{k=1}^{(t-1)/2} \frac{1}{(\zeta_t^k q;q)_\infty (\zeta_t^{-k}q;q)_\infty} 
\mylabel{eq:gHLgenfunc}\\
&= \frac{(q;q)_\infty}
         {(q^t;q^t)_\infty}.
\end{align}
From \eqn{Epent} we have
\beq
\sum_{\vec{\pi}\in\mathcal{P}^{t-1}}  \zeta_t^{\gHLmr} q^{24\abs{\vec{\pi}}+1} 
= 
\frac{1}{(q^{24t};q^{24t})_\infty}
\sum_{n=-\infty}^\infty (-1)^n q^{(6n+1)^2}.
\mylabel{eq:gHLgenfunc2}
\eeq
We see that in the $q$-expansion on the right side of \eqn{gHLgenfunc2} 
the coefficient of $q^n$ is zero when  $n$ is a quadratic nonresidue mod $t$.
Theorem \thm{thm6} (i) follows.

For  part (ii) of Theorem \thm{thm6} we only need Jacobi's triple
product identity \eqn{jtprod}.
We have
\begin{align}
\sum_{\vec{\pi}\in\mathcal{P}^{t-3}}  \zeta_t^{\gHLmr} q^{\abs{\vec{\pi}}}
&= \prod_{k=1}^{(t-3)/2} \frac{1}{(\zeta_t^k q;q)_\infty (\zeta_t^{-k}q;q)_\infty}
\mylabel{eq:gHLgenfunc3}\\
&= \frac{(\zeta_t^{(t-1)/2}q;q)_\infty (\zeta_t^{-(t-1)/2}q;q)_\infty (q;q)_\infty}
         {(q^t;q^t)_\infty}\nonumber\\
&= \frac{\displaystyle \sum_{m=-\infty}^\infty (-1)^{m+1}
         (\zeta_t^{(m+1)(t-1)/2} - \zeta_t^{-m(t-1)/2}) q^{m(m+1)/2}}
         {(1 - \zeta_t^{(t-1)/2}) (q^t;q^t)_\infty},
\nonumber
\end{align}
and
\begin{align}        
&\sum_{\vec{\pi}\in\mathcal{P}^{t-3}}  \zeta_t^{\gHLmr} q^{8\abs{\vec{\pi}}+1}
\mylabel{eq:gHLgenfunc4}\\
&\quad = \frac{1}{(1 - \zeta_t^{(t-1)/2}) (q^{8t};q^{8t})_\infty}
    \sum_{m=-\infty}^\infty (-1)^{m+1}
         \zeta_t^{-m(t-1)/2}
         (\zeta_t^{(2m+1)(t-1)/2} - 1) q^{(2m+1)^2}.
\nonumber
\end{align}
We see that in the $q$-expansion on the right side of \eqn{gHLgenfunc4}
the coefficient of $q^n$ is zero when  $n$ is not a quadratic residue mod $t$,
i.e.\ when $n$ is either a quadratic nonresidue or $n\equiv0\pmod{t}$.
Theorem \thm{thm6} (ii) follows.

\section{Multicranks} \label{sec:mcranks}

In this section we give some extensions of the bicrank to multipartitions
and provide alternative interpretations for some of the congruences given
in Theorem \thm{thm5}. We define two multicranks. These multicranks are defined
in terms of cartesian products of extended partitions and ordinary partitions.
In Section \sect{bicrank}, we defined the set of extended partitions $\mathcal{E}$
and its associated crank and weight function. Recall from Section \sect{ghlbirank}
that $\mathcal{P}$ denotes the set of ordinary partitions, and $\mathcal{P}\subset\mathcal{E}$.
Let $r$ be a positive even integer.
For an extended multipartition
\beq
\vec{\pi}=(\pi_1,\pi_2,\dots,\pi_{r}) \in
\mathcal{E} \times  \cdots  \times
\mathcal{E} \times \mathcal{P} \times \cdots \times \mathcal{P}
=\mathcal{E}^{r/2} \times \mathcal{P}^{r/2},
\eeq
we define multicrank-I by
\beq
\mbox{multicrank-I}(\vec{\pi}) = 
\sum_{k=1}^{r/2} k\cdot\mbox{crank}(\pi_k).
\mylabel{eq:mcrankIdef}
\eeq

For an extended multipartition
\beq
\vec{\pi}=(\pi_1,\pi_2,\dots,\pi_{r}) \in
\mathcal{E} \times \mathcal{E}  
\times \mathcal{P} \times \cdots \times \mathcal{P}
= \mathcal{E} \times \mathcal{E}  \times \mathcal{P}^{r-2},
\eeq
we define multicrank-II by
\beq
\mbox{multicrank-II}(\vec{\pi}) = 
\sum_{k=1}^{2} k\cdot\mbox{crank}(\pi_k) +
\sum_{k=3}^{r} k\,\left(\#(\pi_k) - \#(\pi_{r-k+3})\right).
\mylabel{eq:mcrankIIdef}
\eeq
We note that the $\mbox{bicrank}_2$ corresponds to the multicrank-II when
$r=2$.

For both types of extended multipartitions we define
a sum of parts function and a weight function in the natural way
\beq
\abs{\vec{\pi}} = \sum_{k=1}^{r} \abs{\pi_k},\quad\mbox{and}\qquad
w(\vec{\pi}) = \prod_{k=1}^{r} w(\pi_k).
\mylabel{eq:ewmdef}
\eeq

We have                
\beq
\sum_{\substack{\vec{\pi}\in \mathcal{E}^{r/2}\mathcal{P}^{r/2}\\
      \abs{\vec{\pi}}=n}}  w(\vec{\pi}) 
=
\sum_{\substack{\vec{\pi}\in \mathcal{E}^{2}\mathcal{P}^{r-2}\\
      \abs{\vec{\pi}}=n}}  w(\vec{\pi}) 
= p_{-r}(n).
\mylabel{eq:totw1}
\eeq

\begin{theorem}
\label{thm:thm7}

Let $t>3$ be prime.

\begin{enumerate}
\item[(i)]
The residue of the  multicrank-I mod $t$
divides the extended multipartitions of $n$ from 
$\mathcal{E}^{(t-1)/2}\times\mathcal{P}^{(t-1)/2}$ 
into $t$ equal classes of equal weight 
provided $24n+1$ is a quadratic nonresidue mod $t$.
\item[(ii)]
The residue of the  multicrank-I mod $t$
divides the extended multipartitions of $n$ from 
$\mathcal{E}^{(t-3)/2}\times\mathcal{P}^{(t-3)/2}$ 
into $t$ equal classes of equal weight 
provided   $8n+1$ is not a quadratic residue mod $t$.
\item[(iii)]
The residue of the  multicrank-II mod $t$
divides the extended multipartitions of $n$ from 
$\mathcal{E}^{2}\times\mathcal{P}^{t-5}$ 
into $t$ equal classes of equal weight 
provided   $8n+1$ is a quadratic nonresidue mod $t$.
\end{enumerate}
\end{theorem}

In view of \eqn{totw1}, Theorem \thm{thm7} (i), (ii) provides
alternative combinatorial intepretations of our congruences
for multipartitions given in Theorem \thm{thm5} (i), (ii).
The result in part (iii) is weaker than (ii). We include it
since it is a generalization of the $\mbox{bicrank}_2$.

The proof of Theorem \thm{thm7} is very similar to Theorem \thm{thm6}.
We have
\begin{align}
\sum_{\vec{\pi}\in
\mathcal{E}^{(t-1)/2}\times\mathcal{P}^{(t-1)/2}}
 \zeta_t^{\mcrI} w(\vec{\pi}) q^{\abs{\vec{\pi}}}
&= \left(\prod_{k=1}^{(t-1)/2} F(\zeta_t^k,q) \right) \frac{1}{E^{(t-1)/2}(q)}
\mylabel{eq:mcIgenfunc}\\
&= \prod_{k=1}^{(t-1)/2} \frac{1}{(\zeta_t^k q;q)_\infty (\zeta_t^{-k}q;q)_\infty}
\nonumber\\
&= \frac{(q;q)_\infty}
         {(q^t;q^t)_\infty},
\nonumber
\end{align}
where $F(z,q)$ is the crank generating function given in \eqn{Mvid}.
Theorem \thm{thm7} (i) then follows from \eqn{gHLgenfunc2}. 

Similarly,
\begin{align}
\sum_{\vec{\pi}\in
\mathcal{E}^{(t-3)/2}\times\mathcal{P}^{(t-3)/2}}
 \zeta_t^{\mcrI} w(\vec{\pi}) q^{\abs{\vec{\pi}}}
&= \left(\prod_{k=1}^{(t-3)/2} F(\zeta_t^k,q) \right) \frac{1}{E^{(t-3)/2}(q)}
\mylabel{eq:mcIgenfunc2}\\
&= \prod_{k=1}^{(t-3)/2} \frac{1}{(\zeta_t^k q;q)_\infty (\zeta_t^{-k}q;q)_\infty}.
\nonumber
\end{align}
Theorem \thm{thm7} (ii) then follows from \eqn{gHLgenfunc3}, \eqn{gHLgenfunc4}.

We have
\begin{align}
\sum_{\vec{\pi}\in
\mathcal{E}^{2}\times\mathcal{P}^{t-5}}
 \zeta_t^{\mcrII} w(\vec{\pi}) q^{\abs{\vec{\pi}}}
&= F(\zeta_t,q)\,F(\zeta_t^2,q)\,
 \prod_{k=3}^{(t-1)/2} \frac{1}{(\zeta_t^k q;q)_\infty (\zeta_t^{-k}q;q)_\infty}.
\mylabel{eq:mcIIgenfunc}\\
& = \frac{(q)_\infty^3}{(q^t;q^t)_\infty}.
\nonumber
\end{align}
Then
\beq
\sum_{\vec{\pi}\in
\mathcal{E}^{2}\times\mathcal{P}^{t-5}}
 \zeta_t^{\mcrII} w(\vec{\pi}) q^{8\abs{\vec{\pi}}+1}
= \frac{1}{E(q^{8t})} \sum_{n\ge0}(-1)^n (2n+1) q^{(2n+1)^2},
\mylabel{eq:mcIIgenfunc2}
\eeq
and
Theorem \thm{thm7} (iii) follows.

\section{Concluding Remarks} \label{sec:end}

The two main results of this paper are the combinatorial inpretations
of the 2-colored partition congruences \eqn{2colcongs} in terms
of the Dyson-birank and the $5$-core-birank. The author has been unable to
extended these two results to higher dimensional multipartitions. The extensions
of the Hammond-Lewis birank and Andrews bicrank are much easier because the
generating functions involved are simple infinite products. 

It seems unlikely that a combinatorial proof of \eqn{T5dissect}
is possible. This identity gives the $5$-dissection of the $5$-core-crank
generating function when $z=\zeta_5$. The proof given in the paper relies
on a heavy use of the theory of modular functions. A more elementary
proof is desirable. In \cite{Ga-Ki-St}, a combinatorial proof is given that the 
residue of $5$-core-crank mod $5$ divides the $5$-cores of $5n+4$ into $5$
equal classes. It would interesting to see if the methods of \cite{Ga-Ki-St}
could be extended to give a combinatorial proof of Theorem \thm{thm3},
which is our result for the $5$-core-birank.

It is clear that the generalized-Hammond-Lewis multiranks and our multicranks
are related. For instance, from equations \eqn{gHLgenfunc}, \eqn{mcIgenfunc},
\eqn{gHLgenfunc3}, and \eqn{mcIgenfunc2} we have
\begin{align}
\sum_{\vec{\pi}\in\mathcal{P}^{t-1}}  \zeta_t^{\gHLmr} q^{\abs{\vec{\pi}}}
&=
\sum_{\vec{\pi}\in
\mathcal{E}^{(t-1)/2}\times\mathcal{P}^{(t-1)/2}}
 \zeta_t^{\mcrI} w(\vec{\pi}) q^{\abs{\vec{\pi}}},
\mylabel{eq:gHLmcrid1}\\
\sum_{\vec{\pi}\in\mathcal{P}^{t-3}}  \zeta_t^{\gHLmr} q^{\abs{\vec{\pi}}}
&=
\sum_{\vec{\pi}\in
\mathcal{E}^{(t-3)/2}\times\mathcal{P}^{(t-3)/2}}
 \zeta_t^{\mcrI} w(\vec{\pi}) q^{\abs{\vec{\pi}}}.
\mylabel{eq:gHLmcrid2}
\end{align}
It would interesting to find a combinatorial proof these identities.
However what would be more interesting is to find bijective proofs of
Theorems \thm{thm4} and \thm{thm6}. This is a reasonable problem since
the generating functions involved are simple infinite products.


\noindent
\textbf{Acknowledgement}

\noindent
I would like to thank \dots.


\bibliographystyle{amsplain}

\end{document}